\documentclass[12pt]{amsart}
\usepackage{amsmath}
\usepackage{amssymb}
\usepackage{latexsym}
\def\ds{\displaystyle}
\newenvironment{proof*}{\vskip 2mm\noindent {}}{\hfill $\Box$ \vskip 2mm}
\newtheorem{theorem}{Theorem}
\newtheorem{conjecture}[theorem]{Conjecture}
\newtheorem{corollary}[theorem]{Corollary}
\newtheorem{lemma}[theorem]{Lemma}
\newtheorem{proposition}[theorem]{Proposition}

\newcommand{\C}{{\mathbb{C}}}

\newcommand{\Z}{{\mathbb{Z}}}

\newcommand{\N}{{\mathbb{N}}}
\newcommand{\D}{{\mathbb{D}}}
\renewcommand{\O}{{\mathcal{O}}}
\newcommand{\eps}{\varepsilon}
\newcommand{\Ker}{\operatorname{Ker}}

\title[The Kobayashi--Royden pseudometric of the spectral ball]{On the zero
set of the Kobayashi--Royden pseudometric of the spectral unit ball}
\author{Nikolai Nikolov, Pascal J. Thomas}

\address{Institute of Mathematics and Informatics\\ Bulgarian Academy of
Sciences\\ Acad. G. Bonchev 8, 1113 Sofia,
Bulgaria}\email{nik@math.bas.bg}

\address{Laboratoire Emile Picard, UMR CNRS 5580\\
Universit\'e Paul Sabatier, 118 Route de Narbonne\\ F-31062 Toulouse
Cedex, France} \email{pthomas@cict.fr}



\begin{document}

\begin{thanks}{The initial version of this paper was written during the
stay of the first named author at the Paul Sabatier University,
Toulouse in May--June, 2007.}
\end{thanks}

\begin{abstract} Given $A\in\Omega_n,$ the $n^2$-dimensional
spectral unit ball, we show that if $B$ is in the tangent cone $C_A$
to the isospectral variety at $A,$ then $B$ is a "generalized"
tangent vector at $A$ to an entire curve in $\Omega_n$ if and only
if $B$ is in the tangent cone $C_A$ to the isospectral variety at
$A.$ In the case of $\Omega_3,$ the zero set of the Kobayashi-Royden
pseudometric is completely described.
\end{abstract}

\maketitle

\section{Introduction and results}

Let $\mathcal M_n$ be the set of all $n\times n$ complex matrices.
For $A\in\mathcal M_n$ denote by $sp(A)$ and $\ds
r(A)=\max_{\lambda\in sp(A)}|\lambda|$ the spectrum and the spectral
radius of $A,$ respectively. The spectral ball $\Omega_n$ is the set
$$\Omega_n:=\{A\in\mathcal M_n:r(A)<1\}.$$
The spectral Nevanlinna--Pick problem is the following one: given
$N$ points $a_1,\dots,a_N$ in the unit disk $\Bbb D\subset\Bbb C$
and $N$ matrices $A_1,\dots, A_N\in\Omega_n$ decide whether there is
a mapping $\varphi\in\mathcal O(\Bbb D,\Omega_n)$ such that
$\varphi(a_j )=A_j$, $1\le j \le N$ (cf.
\cite{Agl-You1,Agl-You2,BFT,Cos1,Cos2} and the references there).

The study of the  Nevanlinna--Pick problem in the case $N=2$ reduces
to the computation of the Lempert function, defined as follows for a
domain $D\subset\Bbb C^m$ :
$$ l_D(z,w):=\inf\{|\alpha|:\exists\varphi\in\mathcal O(\Bbb
D,D):\varphi(0)=z,\varphi(\alpha)=w\},\; z,w\in D.$$

The infinitesimal version of the above is the
Carath\'eodory--Fej\'er problem of order $1$: given matrices $A_0,
A_1\in\mathcal M_n$, decide whether there is a mapping $\varphi\in
\mathcal O(\Bbb D,\Omega_n)$ such that $A_0=\varphi(0),$
$A_1=\varphi'(0).$ This problem has been studied in \cite{HMY}.

Its study reduces to the computation of the Kobayashi--Royden
pseudometric, defined as follows for a domain $D\subset\Bbb C^m$:
\begin{align*}
k_D(z;X):=\inf\{|\alpha|:\exists\varphi\in&\mathcal O(\Bbb
D,D):\\&\varphi(0)=z,\alpha\varphi'(0)=X\},\; z\in D,X\in\Bbb C^m.
\end{align*}

To each matrix $A$ we associate its characteristic polynomial
$$P_A(t):= \det(tI-A)=t^n+\sum_{j=1}^n(-1)^j\sigma_j (A)t^{n-j},$$
where $I\in \mathcal M_n$ is the unit matrix,
$$\sigma_j(A):=\sigma_j(\lambda_1,\dots,\lambda_n):=
\sum_{1\le k_1<\dots<k_j\le n}\lambda_{k_1}\dots\lambda_{k_j}$$ and
$\lambda_1,\dots,\lambda_n$ are the eigenvalues of $A.$

Put $\sigma:=(\sigma_1,\dots,\sigma_n):\mathcal M_n\to \Bbb C^n.$
The set
$$\Bbb G_n:=\{\sigma(A):A\in\Omega_n\}$$ is a taut (even hyperconvex)
domain called the symmetrized $n$-disk (cf.
\cite{Agl-You3,Cos2,Jar-Pfl} and references there). Explicit
formulas for $l_{\Bbb G_2}$ and $k_{\Bbb G_2}$ can be found in
\cite{Agl-You2} (see also \cite{Jar-Pfl}) and \cite{HMY},
respectively.

Recall now that a matrix $A\in\mathcal M_n$ is called non-derogatory
if all the blocks in the Jordan form of $A$ have distinct
eigenvalues. Many equivalent properties to this definition may be
found in \cite[Proposition 3]{NTZ}. We point out one of them: $A$ is
non-derogatory if and only if $rank(\sigma_{\ast,A})=n,$ where
$\sigma_{\ast,A}$ stands for the differential of $\sigma$ at the
point $A.$

Denote by $\mathcal C_n$ the open and dense set of all
non-derogatory matrices in $\Omega_n.$

If $A_1,\dots,A_N\in\mathcal C_n$, then any mapping
$\varphi\in\mathcal O(\Bbb D,\Bbb G_n)$ with $\varphi
(\alpha_j)=\sigma(A_j)$ can be lifted to a mapping
$\tilde\varphi\in\mathcal O(\Bbb D,\Omega_n)$ with $\tilde\varphi
(\alpha_j)=A_j,$ $1\le j\le N$ \cite{Agl-You1}. This means that in a
generic case the spectral Nevanlinna--Pick problem for $\Omega_n$
(with dimension $n^2$) can be reduced to the standard
Nevanlinna--Pick problem for $\Bbb G_n$ (with dimension $n$).

By the contractibility of the Lempert function, we have
$$l_{\Omega_n}(A,B)\ge l_{\Bbb G_n}(\sigma(A),\sigma(B)),
\quad A,B\in\Omega_n,$$ and the lifting above implies that equality
holds when both $A,B\in\mathcal C_n.$

By this and the fact that $\mathbb G_n$ is a bounded domain,
$l_{\Omega_n}(A,B)>0$ if $sp(A)\neq sp(B).$ On the other hand, if
$sp(A)=sp(B),$ then there is an entire mapping $\varphi:\Bbb C\to
\Omega_n$ with $\varphi(0)=A$ and $\varphi(1)=B$ \cite{Edi-Zwo};
note that $sp(\varphi(\zeta))=sp(A)$ for all $\zeta \in \mathbb C$,
since by Liouville's Theorem, whenever $\varphi(\C)\subset \Omega,$
then $\sigma \circ \varphi$ is constant. This situation is similar
to that of Brody's theorem for compact manifolds \cite{Bro}: failure
of hyperbolicity (that is, vanishing of the pseudodistance) can be
explained by the presence of a (nonconstant) entire curve in the
manifold.

Restricting again to non-derogatory matrices, a similar lifting
\cite{HMY} implies in the Carath\'eodory-Fej\'er case
$$k_{\Omega_n}(A;B)=k_{\Bbb G_n}(A,\sigma_{\ast,A}(B)),\quad
A\in\mathcal C_n,B\in\mathcal M_n;
$$
in particular, $k_{\Omega_n}(A;B)=0$ if and only if
$\sigma_{\ast,A}(B)=0.$ On the other hand, if
$\sigma_{\ast,A}(B)=0,$ then there is an entire mapping
$\varphi:\Bbb C\to \Omega_n$ with $\varphi(0)=A$ and $\varphi'(0)=B$
(indeed, $B=[Y,A]:=YA-AY$ for some $Y\in\mathcal M_n$ \cite[Proof of
Proposition 3]{NTZ} and then the mapping $\zeta\to e^{\zeta
Y}Ae^{-\zeta Y}$ does the job).

The aim of this paper is to study the zeros of $B\to
k_{\Omega_n}(A;B)$ in the remaining case, where $A$ is a derogatory
matrix, and to relate it to the existence of entire curves tangent
to $B$ at the point $A$ (which is an obvious sufficient condition
for $ k_{\Omega_n}(A;B)=0$).

For $A\in\Omega_n$ denote by $C_A$ the tangent cone (cf. \cite[p.
79]{Chi} for this notion) to the isospectral variety
$$
L_A:=\{C\in\Omega_n:sp(C)=sp(A)\},
$$
that is
$$
C_A:=\{B\in\mathcal M_n:\exists 0<c_j\to 0, C_j\in L_A\mbox{ with }
c_j(C_j-A)\to B\}.$$

Observe that $L_A$ is smooth at $A$ if $A\in\mathcal C_n;$ then
$C_A=\ker\sigma_{\ast,A}$. When $A\not\in \mathcal C_n$, the rank of
$\sigma_{\ast,A}$ is not maximal, so we have $\dim
\ker\sigma_{\ast,A}>n^2 -n$; by \cite[Corollary, p. 83]{Chi}, $C_A$
is an analytic set with  $\dim C_A = \dim L_A = n^2 -n$, so we have
$C_A\subsetneq\ker\sigma_{\ast,A}.$

The following proposition characterizes  the tangent cone $C_A$ as
the set of "generalized" tangent vectors at $A$ to an entire curve
in $\Omega_n$ through $A$ (therefore contained in $L_A$).

\begin{proposition}\label{p1} Let $A\in\Omega_n$ and $B\in\mathcal M_n.$ Then
there are $m\in\Bbb N$\footnote{$\Bbb N:=\{m\in\Bbb Z:m>0\}$}, $m\le
n!$, and $\varphi\in\mathcal O(\Bbb C,\Omega_n)$ with
$\varphi(0)=A,\ \varphi'(0)=\dots=\varphi^{(m-1)}(0)=0,\
\varphi^{(m)}(0)=B$ if and only if $B \in C_A$.
\end{proposition}

Proposition \ref{p1} implies that $C_A$ is contained in the zero set
of the singular Kobayashi pseudometric (cf. \cite{Yu})
\begin{align*}\hat k_{\Omega_n}(A;B)=\inf\{|\alpha|:\exists
m\in\Bbb N,\varphi\in&\mathcal O(\Bbb
D,\Omega_n):\\&\hbox{ord}_0(\varphi-z)\ge m,\
\alpha\varphi^{(m)}(0)=m!X\}.
\end{align*}

A consequence of the proof of Proposition \ref{p1} is the following.

\begin{corollary}\label{4} Let $A\in\Omega_n$ and $B\in C_A.$ Then
the following conditions are equivalent:

(a)  there is $\varphi\in\mathcal O(\Bbb C,\Omega_n)$ with
$\varphi(0)=A$ and $\varphi'(0)=B.$

(b) there are $(r_j)\to\infty$ and $\varphi_j\in\mathcal O(r_j\Bbb
D,\Omega_n),$ $j\in\Bbb N,$ uniformly bounded near $0,$ such that
$\varphi_j(0)=A$ and $\varphi_j'(0)=B.$

(c) there are $r>0$ and $\varphi\in\mathcal O(r\Bbb D,L_A)$ with
$\varphi(0)=A$ and $\varphi'(0)=B.$

\end{corollary}

Before stating the next proposition, we shall define an algebraic
cone $ C'_A \subset \mathcal M_n,$ $A\in\Omega_n.$

For a function $g$ holomorphic near $A$, and $X$ in a neighborhood
of $A$, let $g(X)-g(A) =g_A^*(X-A)+\cdots$, where $g_A^*$ stands for
the homogeneous polynomial of lowest nonzero degree in the expansion
of $g$ near $A$ (the omitted terms are thus of higher order).

Set
$$C^\ast_A:=\{B\in\mathcal M_n:(\sigma_1)_A^*(B)=0,\dots,(\sigma_n)_A^*(B)=0\},$$
$$C'_A:=\cap_{\lambda\in
sp(A)}(\Phi_\lambda)_{\ast,A}^{-1}(C^\ast_{\Phi_{\lambda}(A)}),$$
where
\begin{equation}
\label{autom}\Phi_\lambda(A):= (A-\lambda
I)(I-\overline{\lambda}A)^{-1}.
\end{equation}

Note that
$$C_A\subset C_A^\ast\subset\ker\sigma_{\ast,A}.$$
For the first inclusion, cf. \cite[p. 86, lines 4--6]{Chi}), and for
the second one, use that $$\ker(\sigma_{\ast,A}) =
\{(\sigma_j)_A^*=0 \mbox{ for all }j \mbox{ such that } \deg
(\sigma_j)_A^*=1 \}.$$ Since $C_A$ and $\ker\sigma_{\ast,A}$ are
invariant under automorphisms of $\Omega_n,$ it follows that
$$C_A\subset C_A'\subset\ker\sigma_{\ast,A}.$$

Moreover, if $\dim C_A=\dim  C^*_A,$ that is, $\dim C^*_A =n^2-n,$
then $C_A=C^\ast_A=C_A'$ (cf. \cite[p. 112, Corollary 2]{Chi}).

\begin{proposition}\label{nec}
Let $A\in\Omega_n \setminus \mathcal C_n.$

(i) If $\hat k_{\Omega_n}(A;B)=0$, then $B \in C'_A$;

(ii) $C'_A \neq\ker\sigma_{\ast,A}$.
\end{proposition}

\noindent{\bf Remark.} The cone $C^*_A$ may coincide with
$\ker\sigma_{\ast,A}$ for some $A\in\Omega_n\setminus\mathcal C_n,$
$n\ge 3.$ For example, if $A:=\mbox{diag}(t,\dots,t,0),$ $t\in\Bbb
D_\ast,$ then
$$C^*_A=\ker\sigma_{\ast,A}=\{B\in\mathcal
M_n:tr B=b_{nn}=0\}.$$

The main consequence of Proposition \ref{nec} is the fact that for
$A\in\Omega_n \setminus \mathcal C_n $ and
$B\in\ker\sigma_{\ast,A}\setminus C'_{A},$ a lifting for the
corresponding Carath\'eodory-Fej\'er problem is not possible and
$k_{\Omega_n}(\cdot;B)$ is not a continuous function at $A.$ This
generalizes previous discontinuity results (see \cite{NTZ} and
references therein).

Note also that the cone $ C'_A$ may coincide with $C_A$ in some
cases, for example, for any $A\in\Omega_2$ (then $C_A^\ast=C_A,$
too) and any $A\in\Omega_3$ (for the last, see Proposition \ref{2}
and the discussion before it). We do not know whether this holds in
general. On the other hand, it is not hard to find cases where $C_A
\subsetneq C^*_A$.

\begin{proposition}
\label{conesex} For any $n\ge 3$ there is $A \in\Omega_n$ such that
$C_A \subsetneq  C^*_A$.
\end{proposition}

Now, we state a conjecture about the zero set of $k_{\Omega_n}.$

\begin{conjecture}
\label{mainconj} $k_{\Omega_n}(A;B)=0$ if and only if there is
$\varphi\in \mathcal O(\mathbb C,\Omega_n)$ with $\varphi(0)=A$ and
$\varphi'(0)=B$. In particular, if $k_{\Omega_n}(A;B)=0,$ then $B\in
C_A.$
\end{conjecture}

Conversely however, there are matrices $B\in C_A$ such that
$k_{\Omega_n}(A;B) \neq 0$ (see Proposition \ref{2}(ii) and
Corollary \ref{3}).

There are some cases where our conjecture can be checked.

For example, since $\Omega_n$ is a balanced domain, one has that
$l_{\Omega_n}(0,\cdot)$ and $k_{\Omega_n}(0;\cdot)$ coincide with
the Minkowski function, that is, with the spectral radius. Thus the
zeros of $k_{\Omega_n}(0;\cdot)$ are exactly the zero-spectrum
matrices, and the set of those matrices is a union of complex lines
through $0$.

Also, if $A$ is a scalar matrix, that is, $A=\lambda I,$
$\lambda\in\Bbb C,$ then $B\in C_A,$ if and only if there is
$\varphi\in\mathcal O(\Bbb C,\Omega_n)$ with $\varphi(0)=A$ and
$\varphi'(0)=B$. To see this, use an automorphism of $\Omega_n$ of
the form $(\ref{autom})$ to reduce to the case $A=0.$

Since the derogatory matrices in $\Omega_2$ are exactly the scalar
matrices, we may choose $m=1$ in Proposition \ref{p1} if $n=2,$ and
$C_A$ coincides both with the zeros of $k_{\Omega_2}(A;\cdot)$ and
with the matrices $B=\varphi'(0)$ for some entire curve $\varphi$ in
$\Omega_2$ (on the other hand, $\ker\sigma_{\ast,A}=\{B\in\mathcal M
_2:tr B=0\}$).

Now we shall study the zero set of $k_{\Omega_3}(A;\cdot),$ when $A$
is a non-scalar derogatory matrix. Using first an automorphism of
the form $(\ref{autom})$  and then an automorphism of the form $C\to
D^{-1}CD$, reduces the problem to the following two cases:
$$A=A_t:=\left(\begin{array}{ccc}
0&0&0\\0&0&0\\0&0&t\\
\end{array}\right),\ t\in\Bbb D_\ast,\quad
A=\tilde A:=\left(\begin{array}{ccc}
0&0&0\\0&0&1\\0&0&0\\
\end{array}\right).$$
It is easy to see that
$$C_{A_t}\subset C^\ast_{A_t}=C'_{A_t}=\{B\in\mathcal M_3:
b_{33}=b_{11}+b_{22}=b_{11}^2+b_{12}b_{21}=0\},$$
$$C_{\tilde A}\subset C'_{\tilde A}=\{B\in\mathcal M_3:
b_{11}+b_{22}+b_{33}=b_{32}=b_{12}b_{31}=0\}$$ (to prove, for
example, the second inclusion, use that if $B_\eps=A+\eps
B+o(\eps),$ then $tr B_\eps=\eps tr B+o(\eps),$
$\sigma_2(B_\eps)=-\eps b_{32}+o(\eps)$ and $\det
B_\eps=\eps^2(b_{12}b_{31}-b_{11}b_{32})+o(\eps^2)$). The next
proposition implies, in particular, that
$C_{A_\lambda}=C'_{A_\lambda}$ and $C_{\tilde A}=C'_{\tilde A}$ (use
the fact that the tangent cones are closed or the dimensional
reasoning mentioned above).

\begin{proposition}\label{2} (i) For any $B\in C'_{A_t}$ ($t\neq 0$)
there is $\varphi\in\mathcal O(\Bbb C,\Omega_3)$ with
$\varphi(0)=A_t$ and $\varphi'(0)=B.$

(ii) Let $B\in C'_{\tilde A}.$ Then there is $\varphi\in\mathcal
O(\Bbb C,\Omega_n)$ with $\varphi(0)=\tilde A$ and $\varphi'(0)=B$
if and only if $b_{11}=0$ or $b_{12}\neq b_{31}.$ Otherwise,
$k_{\Omega_3}(\tilde A;B)>0.$
\end{proposition}

\begin{corollary}\label{3} For any $n\ge 3 $ there are $A\in\Omega_n$
and $B\in C_A$ such that $k_{\Omega_n}(A;B)>0.$
\end{corollary}

Since $k_{\Omega_3}(A;B)>0$ if $B\not\in C'_A,$ Proposition \ref{2}
and the discussion before give a complete description of the zero
set of $k_{\Omega_3}.$

\smallskip

Note that the situation is much easier for the
Carath\'eo\-do\-ry-Reiffen pseudometric
$$
\gamma_{\Omega_n}(A;B)=\sup\{|f'(A)B|:f\in\mathcal O(D,\Bbb D)\}.
$$
Here $\gamma_{\Omega_n}(A;B)=0$ if and only if
$\sigma_{\ast,A}(B)=0.$ Indeed, if $\sigma_{\ast,A}(B)\neq 0,$ then
$$\gamma_{\Omega_n}(A;B)\ge\gamma_{\Bbb
G_n}(A;\sigma_{\ast,A}(B))>0.$$ On the other hand, if $A\in\mathcal
C_n$ and $\sigma_{\ast,A}(B)=0,$ then
$$0=k_{\Omega_n}(A;B)\ge\gamma_{\Omega_n}(A;B)\ge 0.$$
It remains to use the density of $\mathcal C_n$ in $\Omega_n$ and
the continuity of the Carath\'eodory-Reiffen pseudometric.
\smallskip

The rest of the paper is organized as follows. The proofs of
Propositions 3 and \ref{conesex} are given in Section 2, the proofs
of Proposition \ref{2} and Corollary \ref{3} -- in Section 3, and
the proofs of Proposition \ref{p1} and Corollary \ref{4} -- in
Section 4.

\medskip

\noindent{\bf Acknowledgment.} The authors wish to thank the referee
for a quick, thorough and perceptive report, thanks to which several
mistakes in the first version have been corrected.

\section{Proofs of Propositions 3 and \ref{conesex}}

\noindent{\bf Proof of Proposition \ref{conesex}.} Set
$A:=\mbox{diag}(0,\dots,0,t,t),$ $t\in \mathbb D_*.$ It is easy to
see that $(\sigma_1)_A^* (B) = \sum_{j=1}^n b_{jj},$ $(\sigma_2)_A^*
(B) = 2t \sum_{j=1}^{n-2} b_{jj} +  t (b_{n-1,n-1} + b_{nn})  ,$
$(\sigma_3)_A^* (B) =t^2 \sum_{j=1}^{n-2} b_{jj}  $. Therefore,
$(\sigma_3)_A^* = t(\sigma_2)_A^* -t^2(\sigma_1)_A^*$ and $\dim
C^\ast_A>n^2-n=\dim C_A$.

\smallskip

\noindent{\bf Proof of Proposition 3.} (i) Let
$\hat\gamma_{\Omega_n}(A;B)$ be the singular Carath\'eo\-dory metric
(cf. \cite{Nik})
$$\hat\gamma_{\Omega_n}(A;B):=\sup\{\left|\frac{f^{(k)}(A)B}{k!}\right|^{1/k}:
k\in\Bbb N,f\in\O(\Omega_n,\D),\hbox{ord}_A f\ge k\},$$ where
$\left|\frac{f^{(k)}(A)B}{k!}\right|=\sum_{|\alpha|=k}D^\alpha
f(A)B^\alpha.$ Since
$$\hat k_{\Omega_n}(A;B)\ge\hat\gamma_{\Omega_n}(A;B),$$ it is enough to
show that $\hat\gamma_{\Omega_n}(A;B)>0$ if $B\not\in C'_A.$ Then
$B\in C^\ast_{\Phi_{\lambda}(A)}$ for some $\lambda\in sp(A).$
Replacing $A$ and $B$ by $\Phi_\lambda(A)$ and
$(\Phi_\lambda)_{\ast,A}(B),$ respectively, we may assume that
$B\not\in C^*_A.$ Then there is $\sigma_j$ such that
$(\sigma_j)_A^*(B)\neq 0.$ Denoting by $k$ the degree of
$(\sigma_j)_A^*,$ it follows that
$$\hat\gamma_{\Omega_n}(A;B)\ge
\left|\frac{(\sigma_j)_A^*(B)}{\binom{n}{j}}\right|^{\frac{1}{k}}>0.$$

(ii) Since $A\in\Omega_n\setminus\mathcal C_n,$ at least two of the
eigenvalues of $A$ are equal, say to $\lambda.$ Applying the
automorphism $\Phi_\lambda$ of $\Omega_n,$ we may assume that
$\lambda=0.$ Since the map $A \to P^{-1} A P$ is a linear
automorphism of $\Omega_n$ for any $P\in\mathcal M_n^{-1}$, we may
also assume that $A$ is in Jordan form. In particular,
$$
A = \left( \begin{array}{cc} A_0 & 0 \\ 0 & A_1 \end{array} \right),
$$
where $A_0 \in \mathcal M_m$, $2\le m \le n$, $sp(A_0)=\{0\}$, $A_1
\in \mathcal M_{n-m}$, $0 \notin sp(A_1)$. Furthermore, there is a
set $J \subsetneq \{2, \dots, m\}$, possibly empty, such that
$a_{j-1,j}=1$ for $j \in J$, and all other coefficients $a_{ij}=0$
for $1\le i, j \le m$. Denote $0\le r := \# J = \mbox{rank} A_0 \le
m-2$.

We set $B := \left( \begin{array}{cc} B_0 & 0 \\ 0 & 0 \end{array}
\right) \in \mathcal M_n,$ where $B_0 = (b_{ij})_{1\le i, j \le m}$
is such that $b_{j-1,j}=-1$ for $j \in \{2, \dots, m\} \setminus J$,
$b_{m1}=1$, and $b_{ij}=0$ otherwise. To complete the proof, it is
enough to show the following.

\begin{lemma}
\label{B} $(\sigma_m)_A^*(B)=1,$ but $\sigma_{\ast,A}(B)=0$.
\end{lemma}

\begin{proof}
We begin by computing $\sigma_j(A_0+hB_0),$ $1\le j \le m,$
$h\in\Bbb C.$ Expanding with respect to the first column, we see
that
$$
\det (tI-(A_0 + h B_0))=t^m+(-1)^{m-1}h^{m-r}.
$$
Comparing the corresponding coefficients of both sides, it follows
that
\begin{equation}\label{h}
\sigma_j(A_0+hB_0)=\left\{\begin{array}{ll}
0,&1\le j\le m-1\\
h^{m-r},&j=m\end{array}.\right.
\end{equation}

Next, we need a general formula for the functions $\sigma_j$. Given
a matrix $M=(m_{ij})_{1\le i, j \le n}$ and a set $E\subset \{1,
\dots, n\}$, we write $\delta_E(M)$ for the determinant of the
matrix $(m_{ij})_{ i, j \in E} \in \mathcal M_{\#E}$. By convention,
$\delta_\emptyset(M)=\sigma_0(M):=1$. Then
\begin{equation}
\label{sigmexp} \sigma_j (M) = \sum_{E\subset \{1, \dots, n\},
\#E=j} \delta_E(M).
\end{equation}

Because of the block structure of our matrices,
$$
\delta_E(A+hB)= \delta_{E\cap \{1, \dots, m\}} (A_0 + h B_0)
\delta_{E\cap \{m+1, \dots, n\}} (A_1).
$$
Therefore
\begin{multline*}
\sigma_j (A+hB) =\sum_{\max(0,j-n+m)\le k \le \min(m,j)} \left(
\sum_{E'\subset \{1, \dots, m\}, \#E'=k} \delta_{E'} (A_0 + h B_0)
\right) \times\\
 \times \left( \sum_{E''\subset \{m+1, \dots, n\}, \#E''=j-k} \delta_{E''} (A_1)
\right)=
\\\sum_{\max(0,j-n+m)\le k \le \min(m,j)}\sigma_k(A_0+hB_0)\sigma_{j-k}(A_1).
\end{multline*}
It follows by $(\ref{h})$ that $\sigma_j(A+hB)= S_1 + S_2$, where
$$S_1=\left\{\begin{array}{ll}
\sigma_j(A_1),&j\le n-m\\
0,&\hbox{otherwise}\end{array},\right.\quad
S_2=\left\{\begin{array}{ll}
h^{m-r} \sigma_{j-m}(A_1),&j\ge m\\
0,&\hbox{otherwise}\end{array}.\right.
$$
In particular,
$$\sigma_j(A)=\left\{\begin{array}{ll}
\sigma_j(A_1),&j\le n-m\\
0,&\hbox{otherwise}\end{array}.\right.$$ Then
$$
\sigma_j (A+hB) - \sigma_j (A) = \left\{\begin{array}{ll}
h^{m-r} \sigma_{j-m}(A_1),&j\ge m\\
0,&\hbox{otherwise}\end{array}\right.
$$
Since $m-r\ge 2$ we get that $\sigma_{\ast,A}(B)=0$, but
$(\sigma_m)_A^*(B)=1.$
\end{proof}

\section{Proofs of Proposition \ref{2} and Corollary \ref{3}}

\noindent{\bf Proof of Proposition \ref{2}.} (i) Let first $B\in
C'_{A_t}.$ We shall write $B$ in the form $B=X+[Y,A_t],$ where $X$
is such that $\psi(\zeta)=A_t+\zeta X\in L_{A_t}$ for any
$\zeta\in\Bbb C.$ Then $\varphi(\zeta)=e^{\zeta
Y}\psi(\zeta)e^{-\zeta Y}$ has the required properties.

It is easy to compute that $\psi(\C)\subset L_{A_t}$ if and only if
$sp(X)=0$ and $x_{11}+x_{22}=x_{11}^2+x_{12}x_{21}=0.$ On the other
hand,
$$[Y,A_t]=t\left(\begin{array}{ccc}
0&0&y_{13}\\0&0&y_{23}\\-y_{31}&-y_{32}&0\\
\end{array}\right).$$

So we may take
$$X=\left(\begin{array}{ccc}
b_{11}&b_{12}&0\\b_{21}&b_{22}&0\\0&0&0\end{array}\right),\quad
Y=t^{-1}\left(\begin{array}{ccc}
0&0&b_{13}\\0&0&b_{23}\\-b_{31}&-b_{32}&0\end{array}\right).$$

\smallskip

(ii) Let first $B\in\C'_{\tilde A}.$ If $b_{11}=0$ or $b_{12}\neq
b_{31},$ it is enough to find (as above) $X$ and $Y$ such that
$B=X+[Y,\tilde A]$ and $\tilde A+\zeta X\in L_{\tilde A}$ for any
$\zeta\in\Bbb C.$ The last condition means that $sp(X)=\{0\}$ and
$x_{32}=x_{12}x_{31}=0.$ On the other hand,
$$[Y,\tilde A]=\left(\begin{array}{ccc}
0&0&y_{12}\\-y_{31}&-y_{33}&y_{22}-y_{33}\\0&0&y_{32}\\
\end{array}\right).$$
Let us assume that $b_{31}=0$ (the computations are similar in the
case $b_{12}=0$). Then we have to choose $X$ of the form
$$X=\left(\begin{array}{ccc} b_{11}&b_{12}&b_{13}-y_{12}
\\b_{21}+y_{31}&b_{22}+y_{32}&b_{23}-y_{22}+y_{33}\\0&0&-b_{11}-b_{22}-y_{32}\\
\end{array}\right)$$
such that $\det X=0$ and $\sigma_2(X)=0,$ that is, $DT=0,$ $D=T^2,$
where we write
$$
D:= \left| \begin{array}{cc} b_{11}&b_{12}
\\b_{21}+y_{31}&b_{22}+y_{32} \end{array} \right|, \quad
T:= b_{11}+b_{22}+y_{32}.
$$
These conditions are satisfied if and only if
$$y_{32}=-b_{11}-b_{22},\quad y_{31}=\left\{\begin{array}{ll}
-b_{21},&b_{11}=0\\
-b_{21}-\frac{b_{11}^2}{b_{12}},&b_{12}\neq 0
\end{array}.\right.$$

It remains to show that if $b_{11}\neq 0$ and $b_{12}=b_{31}=0,$
then $k_{\Omega_3}(\tilde A;B)>0.$ We may assume that $b_{11}=1.$
Set
$$\tilde X=\left(\begin{array}{ccc}
1&0&0\\0&-1&0\\0&0&0\\
\end{array}\right).$$
Choosing $X$ and $Y$ as above, then $B=\tilde X+[Y,A_t].$ Let
$\alpha>0$ and $\varphi\in\mathcal O(\alpha\mathbb D,\Omega_3)$ be
such that $\varphi(0)=A_t$ and $\varphi'(0)=B.$ Setting
$\tilde\varphi(\zeta)=e^{-\zeta Y}\varphi(\zeta)e^{\zeta Y},$ then
$\tilde\varphi\in\mathcal O(\alpha\D,\Omega_3),$
$\tilde\varphi(0)=\tilde A$ and $\tilde \varphi'(0)=X.$ We get that
$k_{\Omega_3}(A_t;B)\ge k_{\Omega_3}(A_t;X).$ The opposite
inequality follows in the same way.

Write $\tilde\varphi$ in the form
$$\tilde\varphi(\zeta)=\tilde A+\zeta\tilde X+\zeta^2\hat X+o(\zeta^2).$$
Then we compute that
$$\sigma_2(\tilde\varphi(\zeta))=\zeta^2(1-\hat x_{32})+o(\zeta^2),\quad
\det\tilde\varphi(\zeta)=-\zeta^3 \hat x_{32}+o(\zeta^3).$$ Since
$|\sigma_2\circ\varphi|<3,$ $|\det\varphi|<1,$ we get by the Cauchy
inequalities that
$$|\hat x_{32}-1|\le 3\alpha^{-2},\quad |\hat x_{32}|\le\alpha^{-3}.$$
So $$k_{\Omega_3}(A_t;B)=k_{\Omega_3}(\tilde A;\tilde
X)\ge\min_{t\in\Bbb C}\max\{\sqrt{|t-1|/3},\root 3\of{|t|}\}>0.$$

\noindent{\bf Proof of Corollary \ref{3}.} Set
$$\tilde A=\left(\begin{array}{ccc}
0&0&0\\0&0&1\\0&0&0\\
\end{array}\right),\
\tilde B_\eps=\left(\begin{array}{ccc}
1&\eps&0\\0&-1&0\\0&0&0\\
\end{array}\right),$$
$$A=\left(\begin{array}{cc}\tilde A&O\\O&O
\end{array}\right),\
B_\eps=\left(\begin{array}{cc}\tilde B_\eps&O\\O&O
\end{array}\right).$$
It follows as in the proof of Proposition \ref{2} (ii) that

$\bullet$ $k_{\Omega_n}(A;B_0)>0;$

$\bullet$ for $\eps\neq 0,$ there is $\varphi_\eps\in\mathcal O(\Bbb
C,\Omega_n)$ with $\varphi_\eps(0)=A$ and $\varphi_\eps'(0)=B_\eps.$

\noindent Then $B_\eps\in C_A,$ $\eps\neq 0,$ and hence $B_0\in
C_A.$

\section{Proofs of Proposition \ref{p1} and Corollary \ref{4}}

\noindent{\bf Proof of Corollary \ref{4}.} The implication
$(a)\Rightarrow (b)$ is trivial and the main implication
$(c)\Rightarrow (a)$ is a partial case of Proposition \ref{entire}
below.

It remains to prove that $(b)\Rightarrow (c).$ Let
$\psi(\zeta)=\sum_{k=0}^\infty A_k\zeta^k\in\mathcal O(s\Bbb
D,\Omega_n)$ for some $s>0.$ Let $a_{k,\psi}\in\Bbb C^{(k+1)\times
n^2}$ be the vector with components the entries of $A_0,\dots,A_k$
(taken in some order). Note that
$$\sigma_l(\psi(\zeta))=\sum_{k=0}^\infty p_{l,k}(a_{k,\psi})\zeta^k,\quad
1\le l\le n,$$ where $p_{j,k}$ are polynomials.

Let now $(r_j)\to\infty$ and $\varphi_j\in\mathcal O(r_j\Bbb
D,\Omega_n),$ $j\in\Bbb N,$ uniformly bounded near $0,$ be such that
$\varphi_j(0)=A$ and $\varphi_j'(0)=B.$ Then we may assume that
$\varphi_j\to\varphi\in\mathcal O(r\Bbb D,\Omega_n)$ for some $r>0.$
Hence $p_{l,k}(a_{k,\varphi_j})\to p_{l,k}(a_{k,\varphi}).$ On the
other hand, $|p_{l,k}(a_{k,\varphi_j})|\le\binom{n}{l}/r_j^k\to 0,$
$k>0,$ by the Cauchy inequalities. Hence $p_{l,k}(a_{k,\varphi})=0,$
$k>0,$ that is,
$\sigma_l(\varphi(\zeta))=\sigma_l(\varphi(0))=\sigma_l(A).$ This
means that $\varphi(\zeta)\in L_A.$
\smallskip

\noindent{\bf Proof of Proposition \ref{p1}.} It is clear that if
such a $\varphi$ exists, then $B\in C_A.$

Conversely, let $B\in C_A$. Then, by \cite[p. 86, Proposition
1]{Chi}, there exists a one-dimensional irreducible analytic variety
$L_{A,B} \subset L_A$, tangent to $B$ at $A$. Now, by \cite[p. 80,
Proposition]{Chi}, there exist $m\in\N$, $r>0$ and
$\psi\in\O(r\D,L_{A,B})$ such that $\psi(\zeta) = A+\zeta^m
B+o(\zeta^{m}).$

The integer $m$ is the number of sheets in the (local) branched
covering provided by the orthogonal projection from $L_A$ to a
suitable linear subspace of dimension $n^2-n$ \cite{Chi}. This
number of sheets corresponds to the cardinality of the solution set,
in each generic fiber of the projection, of the equations
$\sigma_j(M)=\sigma_j(A)$, $1\le j \le n$. Bezout's theorem shows
that this cardinality is less than or equal to the product of the
degrees of the polynomials, so here $m \le n!$.

The above considerations will prove the estimate in Proposition
\ref{p1}, provided that we can replace $\psi$ by an entire map
$\tilde \psi$ with the same expansion up to order $m$ near $0$. So
the proof of Proposition \ref{p1} reduces to the following.

\begin{proposition}
\label{entire} If $A\in\Omega_n,$ $m\in\N$ and $\psi\in\mathcal
O(r\D,L_A)$ for some $r>0,$ then there is $\tilde\psi\in\mathcal
O(\C,L_A)$ with $\tilde\psi(\zeta)=\psi (\zeta)+o(\zeta^m)$.
\end{proposition}

\noindent{\bf Proof.} We want to reduce that problem by replacing
each matrix $\psi(\zeta)$ by a conjugate matrix $\varphi(\zeta)$ (in
particular, they will have the same spectrum, so we remain inside
$L_A$ and inside $\Omega_n$). If we can manage this so that
$\varphi(\zeta)$ is upper triangular, then an entire map with the
same spectrum matching $\varphi$ up to order $m$ can be obtained by
taking the Taylor polynomial of degree $m$ of each coefficient of
$\varphi$.

To proceed with this program, first we need to show that conjugation
(with a holomorphic change of basis) does not change the problem.

Let $\mathcal M_n^{-1}$ stand for the group of all invertible
$n\times n$ matrices.

\begin{lemma}
\label{change} Let $r>0,$ $P \in\mathcal O(r\D, \mathcal M_n^{-1})$
and $\psi\in\mathcal O(r\D,\Omega_n)$.

Write $\varphi(\zeta) := P(\zeta)^{-1}\psi(\zeta)P(\zeta)$, and
assume that there exists $\tilde \varphi\in\mathcal O(\C,\Omega_n)$
such that near $0$, $\tilde
\varphi(\zeta)=\varphi(\zeta)+o(\zeta^m).$

Then there exists $\tilde\psi\in\mathcal O(\C,\Omega_n)$, conjugate
to  $\tilde \varphi$ (in particular they have the same spectrum)
such that near $0$,
$$\tilde\psi(\zeta)=\psi(\zeta)+o(\zeta^m).$$
\end{lemma}

Note once again that Liouville's theorem implies that the entire
maps $\tilde \varphi, \tilde\psi$ actually map to $L_{\tilde
\varphi(0)}=L_{\tilde \psi(0)}$.

\begin{proof} Note first that, because the exponential is locally
Lipschitz, $\exp(A+M)= \exp A + O(M).$

Denote by $L_m(x)$ the Taylor polynomial of degree $m$ at $0$ for
the function $x \to \ln(1+x)$. Since $\exp (\ln(1+x))=1+x$ and
$\ln(1+x)=L_m(x)+o(x^m)$, we have $\exp (L_m(x))=1+x+o(x^m)$. So
$\exp (L_m(A))=I+A+o(A^m)$.

Now write
$$
P(\zeta) = P(0) \left( I + \sum_{k=1}^m A_k \zeta^k + O(\zeta^{m+1})
\right) =: P(0) \left( I + M(\zeta) \right) .
$$
 Define $P_1$ to be the unique matrix-valued polynomial of degree $\le m$
in $\zeta$ so that
$$
L_m \left( \sum_{k=1}^m A_k \zeta^k \right) = P_1 (\zeta)  +
o(\zeta^m) .
$$
Then, remarking that $M(\zeta)=o(1)$, we have
$$\exp (P_1(\zeta)) = \exp (L_m(M(\zeta))+ o(\zeta^m) =
I + M(\zeta)+ o(\zeta^m),$$
 so that $P(0) \exp (P_1(\zeta)) = P(\zeta) + o(\zeta^m).$  Then it is easy to
see that $\exp (-P_1(\zeta)) P(0)^{-1}= P(\zeta)^{-1} + o(\zeta^m),$
and $\tilde P(\zeta) := P(0) \exp (P_1(\zeta))$ defines an entire
map. So $\tilde \psi := \tilde P \tilde\varphi \tilde P^{-1}$
satisfies the requirements.
\end{proof}

We now reduce the proof of Proposition \ref{entire} to the case of
nilpotent matrices (that is, $sp(A)=\{0\}$).

\begin{lemma}
\label{nilpotent} Suppose that Proposition \ref{entire} holds with
the additional hypothesis that $sp(A)=\{0\}$. Then it holds for an
arbitrary matrix $A$.
\end{lemma}

\begin{proof}
Write $sp(\psi(\zeta))=sp(A)=\{\mu_j:1\le j \le k\}$ where the
$\mu_j$ are distinct eigenvalues with respective algebraic
multiplicity $m_j$. Let $S_j(\zeta )=\ker (\psi(\zeta )-\mu_j
I)^{m_j}$ be the associated generalized eigenspace. Choose
$\{e_1,\dots, e_n\}$ a basis of $\C^n$ such that $S_j
(0)=span\{e_i:1\le i-\sum_{l=1}^{j-1}m_l\le m_j\}.$ By Lemma
\ref{change}, without loss of generality we may assume that the
matrices are written in this basis, and therefore $\psi(0)$ is a
block matrix.

By continuity of the various determinants involved, there exists
some $r'>0$ such that for $|\zeta |<r'\le r$, we still
have, for
each $j$,  $\C^n= S_j(\zeta) \oplus \bigoplus_{l:l\neq j} S_l(0)$.
Then there is a unique linear projection $\pi_{j,\zeta}$ defined on  $\C^n$
such that $\pi_{j,\zeta} (\C^n) = S_j(\zeta)$ and
$\Ker \pi_{j,\zeta} = \bigoplus_{l:l\neq j} S_l(0)$.  This restricts
to a linear isomorphism from $S_j(0)$ to $S_j(\zeta)$. Therefore the vectors $\{ \pi_{j,\zeta} (e_i), 1\le i-\sum_{l=1}^{j-1}m_l\le m_j \}$,
being obtained as solution of a Cramer system of linear equations with holomorphic coefficients, depend holomorphically on $\zeta$ in
$D(0,r')$. Thus $\{ \pi_{j,\zeta} (e_i), 1\le i-\sum_{l=1}^{j-1}m_l\le m_j, 1 \le j \le k \}$ form a basis of $\C^n$ adapted to the direct sum
decomposition in the $S_j(\zeta)$. If we write $P(\zeta)$ for the matrix of the coordinates of the vectors of this new basis
expressed in the standard basis, it depends holomorphically on $\zeta$ in
$D(0,r')$
and the new matrix $\hat \psi (\zeta ):=
P(\zeta)^{-1}\psi (\zeta )P(\zeta)$, has the same block structure
as $\psi(0)$:
$$
\hat \psi (\zeta )= \left(
\begin{array}{cccc}
\hat \psi_1 (\zeta ) & 0 & \cdots & 0 \\
0 & \hat \psi_2 (\zeta ) & \cdots & 0 \\
\vdots & \vdots & \ddots & \vdots \\
0 & 0 & \cdots & \hat \psi_k (\zeta )
\end{array}\right) ,
$$
where $\hat \psi_j \in \mathcal O (r' \D, \Omega_{m_j})$, $sp(\hat
\psi_j (\zeta )) = \{\mu_j\}$. The map $\omega_j$ defined by
$$
\omega_j(\zeta ): = \left( \mu_j I_{m_j} - \hat \psi_j (\zeta )
\right) \left( I_{m_j} - \bar \mu_j \hat \psi_j (\zeta )\right)^{-1}
$$
is in $\mathcal O (r' \D, \Omega_{m_j})$, and its values are
nilpotent matrices. By our hypothesis there are maps $\tilde
\omega_j \in \mathcal O (\C, \Omega_{m_j})$ such that $\tilde
\omega_j (\zeta)=  \omega_j (\zeta) +o(\zeta^m)$ (and therefore with
nilpotent values). Define
$$
\tilde \psi_j (\zeta ): = \left( \mu_j I_{m_j} -\tilde \omega_j
(\zeta) \right) \left(  I_{m_j} - \bar \mu_j \tilde \omega_j (\zeta)
\right)^{-1},
$$
and
$$
\tilde \psi (\zeta ):= \left(
\begin{array}{cccc}
\tilde \psi_1 (\zeta ) & 0 & \cdots & 0 \\
0 & \tilde \psi_2 (\zeta ) & \cdots & 0 \\
\vdots & \vdots & \ddots & \vdots \\
0 & 0 & \cdots & \tilde \psi_k (\zeta )
\end{array}\right) .
$$
It is easy to see that $\tilde \psi \in  \mathcal O (\C,
\Omega_{n})$ and that $\tilde \psi (\zeta ) = \hat \psi (\zeta )
+o(\zeta^m)$.
\end{proof}

\begin{lemma}
\label{triangle} If $m\in\N$ and $\psi\in\mathcal O(r\D,L_0)$ for
some $r>0,$ then there is $P \in\mathcal O(r'\D, \mathcal M_n^{-1})$
such that $\varphi(\zeta) := P(\zeta)^{-1}\psi(\zeta)P(\zeta)$ is a
strictly upper triangular matrix for all $\zeta \in r'\D$.
\end{lemma}

Proposition \ref{entire} follows from Lemma \ref{triangle}.

Indeed, by Lemma \ref{nilpotent}, we can make the
additional hypothesis that $\psi$ has nilpotent values, that is,
$\psi\in\mathcal O(r\D,L_0)$. By Lemma \ref{triangle},
there is, for every $\zeta$ in a neighborhood of $0$,
a strict upper triangular matrix $\varphi (\zeta)
= P(\zeta)^{-1}\psi(\zeta)P(\zeta)$. If we
replace each of the holomorphic coefficients $\varphi_{ij}(\zeta )$,
$1\le i<j\le n,$ by its Taylor polynomial of order $m$,
$\tilde \varphi_{ij}(\zeta ) := \sum_{k=0}^m \varphi_{ij}^{(k)}(0) \zeta^k/k!
$, we obtain $\tilde\varphi$ an approximation up to order $m$ of our mapping which is
entire, and still strictly upper-triangular, therefore still with
spectrum reduced to $0$. Since $P(\zeta)$ depends holomorphically
on $\zeta$ in a neighborhood of $0$, we may apply
 Lemma \ref{change} and obtain the matrix $\tilde \psi(\zeta)$
still with spectrum reduced to $0$, approximating $\psi$ to order $m$, and entire in $\zeta$.

\noindent{\bf Proof of Lemma \ref{triangle}.} We are working with a
matrix $\psi(\zeta )$ which satisfies $\psi(\zeta )^n=0$ for all
$\zeta \in r\Bbb D$. For $1\le k \le n$, let $r_k(\zeta ):=
\mbox{rank} \left( \psi(\zeta )^k \right)$. For a matrix $M$, $
\mbox{rank} (M)\le l<n$ if and only if all the minor determinants of
size $l+1$ vanish; in our case those determinants are holomorphic
functions of $\zeta $, therefore
$$r_k(\zeta)=\max_{\D}r_k=:\tilde r_k$$
for all $\zeta\in r\D$ except on a discrete set. By replacing $r$ by
a smaller positive number, if needed, we may assume that $r_k(\zeta
)=\tilde r_k$ for all $\zeta\in r\D_\ast:=r\D\setminus\{0\}$. Set
$$n_k:=n-\tilde r_k =\dim\ker\psi(\zeta)^k,\quad
\zeta\in r\D_\ast.$$ It is a classical fact from linear algebra
that, for a nilpotent matrix, if $p:=\min\{k: n_k=n\}$, then $1\le
n_1<\cdots<n_p=n_{p+1}=\cdots=n$.

For $1\le k \le p$ and $\zeta\in r\D_\ast,$ set $V_k (\zeta ):=
\ker\psi(\zeta )^k$. Since the Grassmannian $\mathcal G(n,n_k)$ is
compact, we may find a sequence $(\zeta _i)\to 0$ and vector
subspaces $V_k(0) \in \mathcal G(n,n_k)$ such that
$\lim_{i\to\infty} V_k(\zeta_i)=V_k(0)\subset\ker\psi(0)^k $, $1\le
k \le p$.

Our problem will be solved if we find $\eps >0$ and holomorphic
mappings $v_j \in\mathcal O(\eps\D,\C^n)$ such that $\{v_j(\zeta ),
1 \le j \le n_k\}$ is a basis of $V_k(\zeta )$, $\zeta\in\eps\D,$
$1\le k \le p$. (In particular, $\lim_{\zeta \to 0} V_k(\zeta
)=V_k(0)$).

We shall proceed by induction on $k$, and on $j$ for each fixed $k$.
The value of $\eps$ may be reduced at each step, but we keep the
same notation.

By convention we will set $n_0=0$, and consider $\emptyset$ as a
basis of $\{0\}$. Suppose that we already have determined $\{v_j:1 \le
j\le n_k\}$. Choose an $r_{k+1}\times r_{k+1}$ minor matrix of
$\psi^{k+1}$ whose determinant, denoted $\delta_{k+1}$,
is holomorphic and does not
vanish on $\eps\D_\ast$ and eliminate the unknowns corresponding
to the columns of this minor; the other unknowns are then expressed
in terms of the former with coefficients which are rational in the
coefficients of the matrix $\psi^{k+1}$, so that we obtain
meromorphic
vector-valued functions $u_i$ on $\eps\D$ so that $\{u_i(\zeta), \le
j\le n_{k+1}\}$ is a basis of
$V_{k+1}(\zeta )$ for $\zeta\in\eps\D_\ast.$ Those
functions are of the form $u_i:=f_i/\delta_{k+1}$,
where $f_i\in \mathcal O(\eps\D,\C^n)$.

By linear algebra, for each fixed $\zeta$, there exist
a set $I(\zeta )$ of $n_{k+1}-n_k$ indices $i$ so that
$\{ v_j(\zeta ),u_i(\zeta):1 \le j\le
n_k, i \in I(\zeta )\}$ form a basis of $V_{k+1}(\zeta )$.

Using the fact that all determinants that we have to compute
to determine the rank of a system of vectors are
meromorphic in $\eps\D$ and reducing $\eps$ if needed, we may choose
a fixed set $I$ so that
$\{ v_j(\zeta ),u_i(\zeta):1 \le j\le
n_k, i \in I\}$ form a basis of $V_{k+1}(\zeta )$, for all $\zeta\in\eps\D_\ast.$
Re-index the functions
$\{u_i, i \in I\}$ as $w_j, $ $n_k+1 \le j
\le n_{k+1}$: we have that $\{ v_j(\zeta ),w_l(\zeta):1 \le j\le
n_k<l\le n_{k+1}\}$ form a basis of $V_{k+1}(\zeta )$ for
$\zeta\in\eps\D_\ast.$

Since $\delta_{k+1}$ vanishes at most at $0$, we can multiply each
vector-valued function  $w_l$ by $\zeta^{\alpha_l}$, with $\alpha_l\in\Z$ chosen
such that $\zeta^{\alpha_l}w_l(\zeta)$ extends to a map $\tilde
w_l\in\mathcal O(\eps\D,\C^n)$ with $\tilde w_l(0)\neq 0$.

We need to modify the $\tilde w_l$ to ensure that we still have a
system of maximal rank at the origin (this is in the spirit of the
Gram-Schmidt orthogonalization process). We proceed by induction on
$j\ge n_k+1.$ Suppose we have determined $v_j$ as in the statement
of the lemma up to some $j_0\ge n_k$ such that
\begin{multline*}
span\{v_j(\zeta ):1 \le j \le j_0\}
\\
=span\{ v_j(\zeta ),\tilde w_l(\zeta):1\le j\le n_k<l\le
j_0\},\quad\zeta\in\eps\D_\ast,
\end{multline*}
$$rank\{ v_j(\zeta ):1 \le j \le j_0\} = j_0,\quad
\zeta\in\eps\D.$$ Let $W(\zeta ):=span\{v_1(\zeta
),\dots,v_{j_0}(\zeta),\tilde w_{j_0+1}(\zeta )\},$ $\zeta\in
\eps\D_\ast.$ Then $\dim W(\zeta )\\=j_0+1$. Again by using the
compactness for the Grassmannian, we may choose a sequence $(\zeta
_i)\to 0$ such that $\lim_{i\to\infty}W(\zeta _i)=: W(0)$ exists.
Since all the $v_j$ and $\tilde w_{j_0+1}$ are continuous at $0$, we
easily deduce that $span\{v_1(0),\dots,v_{j_0}(0),\tilde
w_{j_0+1}(0)\}\subset W(0)$. Choose a vector $w$ such that
$\{v_1(0),\dots, v_{j_0}(0),w\}$ form a basis of $W(0)$.

By reducing $\eps$ if needed, we may assume that all the space
$W(\zeta )$ and $W(0)$ are each in direct sum with a common
supplementary space $Y$. Then define $v_{j_0+1}( \zeta )$ to be the
projection of the fixed vector $w$ to $W(\zeta )$, parallel to $Y$,
for $\zeta \neq 0$; in particular  $v_{j_0+1}(0)=w.$ Since
$v_{j_0+1}( \zeta )$ is obtained by solving a system of linear
equations with unique solution, it is meromorphic in a neighborhood
of $0$, but it is also clearly continuous near $0$, hence
holomorphic. Since the system $\{v_1(0),\dots,
v_{j_0}(0),v_{j_0+1}(0)\}$ is independent, then $\{v_1(\zeta),\dots,
v_{j_0}(\zeta),v_{j_0+1}(\zeta)\}$ also is independent for $\zeta$
small enough, and since it is contained in $W(\zeta )$ by
construction, it forms a basis of that subspace. Reducing $\eps$ yet
again if needed, the induction may proceed.

\end{document}